\documentclass[a4paper,12pt]{amsart}

\usepackage{t1enc}

\def\Bbb{\mathbb}
\def\Cal{\mathcal}

\newcommand{\wh}{\widehat}

\newcommand{\cq}{{\Cal Q}}

\newcommand{\K}{{\Bbb K}}

\newcommand{\cL}{{\Cal L}}

\newcommand{\Ric}{\operatorname{Ric}}

\newcommand{\om}{\omega}
\newcommand{\Om}{\Omega}

\newcommand{\newc}{\newcommand}

\newcommand{\si}{\sigma}

\newcommand{\De}{\Delta}

\newcommand{\Up}{\Upsilon}

\newc{\al}{\mbox{\boldmath$ \Delta$}}

\newtheorem{theorem}{Theorem}[section]
\newtheorem{lemma}[theorem]{Lemma}
\newtheorem{proposition}[theorem]{Proposition}

\newcommand{\ce}{{\Cal E}}

\usepackage{amssymb}
\usepackage{amscd}

\newcommand{\nd}{\nabla}

\newcommand{\Rho}{P}

\newcommand{\Pa}{{\Bbb I}}

\newcommand{\nn}[1]{(\ref{#1})}

\newcommand{\D}{\mbox{\boldmath{$ D$}}}
\newcommand{\tU}{\tilde{U}}

\newcommand{\bg}{\mbox{\boldmath{$ g$}}}

\newcommand{\V}{P}
\newcommand{\J}{J}

\newc{\aR}{\mbox{\boldmath{$ R$}}}
\newc{\aS}{\mbox{\boldmath{$ S$}}}
\newc{\aDeR}{\mbox{\boldmath{$ U$}}_B{}^P{}_C{}^Q}
\newc{\aDe}{\mbox{\boldmath$ \Delta$}}
\newc{\aNd}{\mbox{\boldmath$ \nabla$}}

\newc{\aK}{\mbox{\boldmath{$ K$}}}
\newc{\aL}{\mbox{\boldmath{$ L$}}}
\newcommand{\bD}{\mbox{\boldmath{$ D$}}}
\newcommand{\X}{\mbox{\boldmath{$ X$}}}
\newcommand{\h}{\mbox{\boldmath{$ h$}}}

\def\endrk{\hbox{$|\!\!|\!\!|\!\!|\!\!|\!\!|\!\!|$}}

\def\sideremark#1{\ifvmode\leavevmode\fi\vadjust{\vbox to0pt{\vss
 \hbox to 0pt{\hskip\hsize\hskip1em
 \vbox{\hsize3cm\tiny\raggedright\pretolerance10000
 \noindent #1\hfill}\hss}\vbox to8pt{\vfil}\vss}}}%
                        
\author{A. Rod Gover} \email{gover@math.auckland.ac.nz} \title{Laplacian
operators and Q-curvature on conformally Einstein manifolds}
\begin{document}

\begin{abstract}
 A new definition
  of canonical conformal differential operators $P_k$
  ($k=1,2,\cdots)$, with leading term a $k^{\rm th}$ power of the
  Laplacian, is given for conformally Einstein manifolds of any
  signature.  These act between density bundles and, more generally,
  between weighted tractor bundles of any rank. By construction these
  factor into a power of a fundamental Laplacian associated to
  Einstein metrics. There are natural conformal Laplacian operators on
  density bundles due to Graham-Jenne-Mason-Sparling (GJMS).  It is
  shown that on conformally Einstein manifolds these agree with the
  $P_k$ operators and hence on Einstein manifolds the GJMS operators
  factor into a product of second order Laplacian type operators.  In
  even dimension $n$ the GJMS operators are defined only for $1\leq
  k\leq n/2$ and so, on conformally Einstein manifolds, the $P_{k}$
  give an extension of this family of operators to operators of all
  even orders. For $n$ even and $k>n/2$ the operators $P_k$ are each given
  by a natural formula in terms of an Einstein metric but they are not
  natural conformally invariant operators in the usual sense. They are
  shown to be nevertheless canonical objects on conformally Einstein
  structures.  There are generalisations of these results to operators
  between weighted tractor bundles. It is shown that on Einstein
  manifolds the Branson Q-curvature is constant and an explicit
  formula for the constant is given in terms of the scalar curvature.
  As part of development, conformally invariant tractor equations
  equivalent to the conformal Killing equation are presented.
\end{abstract}

\maketitle

\pagestyle{myheadings}
\markboth{A. Rod Gover}{Conformal Laplacian
operators and Q on Einstein manifolds}

\section{Introduction}

Conformally invariant operators with leading term a power of the
Laplacian $\Delta$ have been a focus in mathematics and physics for
over 100 years.  The earliest known of these is the conformally
invariant wave operator which was first constructed for the study of
massless fields on curved space-time (see
e.g. \cite{dirac}).   Its Riemannian signature
elliptic variant, often called the Yamabe operator, controls the
transformation of the Ricci scalar curvature under conformal rescaling
and so plays a critical role in the Yamabe problem on compact
Riemannian manifolds.  A conformal operator with principal part $
\Delta^2$ is due to Paneitz \cite{Pan} (see also \cite{Rie,ESlo}), and
 sixth-order analogues were constructed in
\cite{Br85,Wuensch86}. Graham, Jenne, Mason and Sparling (GJMS) solved
a major existence problem in \cite{GJMS} where they used a formal
geometric construction (the Fefferman-Graham ambient metric) to show
the existence of conformally invariant differential operators $
\Box_{k}$ (to be referred to as the GJMS operators) with principal
part $ \Delta^k$.  In odd dimensions, $ k$ is any positive integer,
while in dimension $ n$ even, $k $ is a positive integer no more than
$ n/2$.  The $ k=1$ and $ k=2$ cases recover, respectively, the Yamabe
and Paneitz operators. Recently Gover and Hirachi have shown that the
GJMS result is sharp in the sense that on general conformal manifolds
there are no other conformal operators between density bundles with
leading term a power of the Laplacian \cite{GoHirnon}.

Intimately linked to this family of operators is Branson's
Q-curvature.  \cite{BrO,Brsharp}. In particular the conformal
variation of Q is by the dimension order GJMS operator acting linearly on the 
conformal factor  and this
leads to a generalisation to higher even dimensions of the scalar
curvature prescription problem in dimension 2. 
This has motivated much recent activity, see for example \cite{CGY,CQY,DM}.

Recently Graham and Hirachi established \cite{GrH} that the
total metric variation of $\int Q$ is the conformally invariant
{\em obstruction} tensor that arises as an obstruction \cite{FGast} to
the Fefferman-Graham ambient metric construction. Since the
obstruction tensor vanishes on conformally Einstein manifolds
\cite{GrH,GoPetobstrn} it follows that conformally Einstein
metrics are critical for $\int Q$.  One of the main results here is a
simple proof of the following theorem which shows that for the Q-curvature
Einstein metrics have an elevated status from within a conformal class.
\begin{theorem}\label{Qthm} On a manifold of even dimension $n$ and with an Einstein metric 
  $g$, the Q-curvature 
is constant and given by
\begin{equation}\label{Qformula}
Q^g=(-1)^{n/2}(n-1)! \Big(\frac{{\rm Sc^g}}{n(n-1)} \Big)^{n/2},
\end{equation}
where ${\rm Sc^g}$ is the Ricci scalar curvature.
\end{theorem} 
\noindent The proof uses another central result of the article, as follows.
\begin{theorem}\label{prodthm} On an $n$-manifold with Einstein metric 
$g$  the order $2k$ GJMS operator is given by
\begin{equation}\label{prodformula}
\Box_{k} = \prod_{l=1}^k(\De^g - c_l {\rm Sc}^g),
\end{equation}
where $c_l= (n+2l-2)(n-2l)/(4n(n-1))$ and $\Delta^g=\nd^a\nd_a$ is the
Laplacian for $g$.
\end{theorem}
\noindent Note that the left (i.e. $l=1$) factor is the conformal Laplacian: 
$Y=\De^g - \frac{(n-2)}{4(n-1)} {\rm Sc}^g$. Observe also that for $k< n/2$ 
we have $c_l>0$ ($1\leq l\leq k$) and (for $n$ even) $c_{n/2}=0$. 

There are obvious spectral applications for Theorem \ref{prodthm}, and
its generalisations discussed below. For example it follows that, in
the Einstein scale $g$, each eigenfunction for $\De^g$, with
eigenvalue $\lambda$ say, is an eigenfunction for the GJMS operator
$\Box_k$ and the corresponding eigenvalue is then given explicitly in
terms of ${\rm Sc}^g$ and $\lambda$ via the formula
\nn{prodformula}. In particular, in the compact Riemannian setting the
eigenfunctions of $\De^g$ form a complete spectral resolution and so the
Theorem gives a complete spectral resolution of each of the GJMS operators
in terms of this. We may alternatively rewrite the product
\nn{prodformula} in terms of the conformal Laplacian $Y$ and similar
comments apply with $Y$ replacing $\De^g$. 

In the cases where $k$ is not half the dimension then
$Q_k:=\frac{2}{n-2k}\Box_k 1$ is sometimes known as the {\em
non-critical Q-curvature} and has a role of yielding curvature
prescription problems that are higher order generalisations of the
classical Yamabe problem \cite{Br85,Brsharp}. From Theorem \ref{prodthm} this is constant
on Einstein manifolds and given by $\frac{2}{n-2k}\prod_{l=1}^k(
- c_l {\rm Sc}^g) $.

Theorem \ref{prodthm} generalises in two ways. These are both captured
in theorem \ref{bigprodthm}, at least given Theorem \ref{natthm}. The
first direction is concerned with the domain and range bundles.  The
GJMS operators are defined in \cite{GJMS} as operators between
density bundles. In \cite{GoPetobstrn} it is shown that these
operators generalise to conformally invariant Laplacian power type
operators between weighted tractor bundles. A class of these is
discussed in Theorem \ref{powerslap} below and we show that these
operators also factor in a way formally identical to the GJMS
operators in Theorem \ref{prodthm} above.  Tractor bundles are
reducible but indecomposable bundles with a natural conformally
invariant connection and are intimately linked to tensor bundles
through the connection soldering form. (The conformal tractor
connection is due to Cartan \cite{Cartan} and, independently, Thomas
\cite{T}. It was rediscovered and put into a modern in framework in
\cite{BEGo}. A complete treatment and generalisations are given in
\cite{CapGotrans}.) The point of the tractor Laplacian operators
developed in \cite{GoPetobstrn} is that they enable the construction
of large classes of other conformally invariant operators through the
curved translation principle of Eastwood et al.\ \cite{ER,Esrni}.
Examples of such translation constructions are treated in for example
\cite{ER,GoPetLap,GoPetobstrn}. There is evidence \cite{GoSilth} that
that every non-standard operator predicted in \cite{ESlo} can be
obtained via translation from the Laplacian tractors in
\cite{GoPetobstrn}. In any case where operators arise via translation
from the operators of Theorem \ref{powerslap}, then Theorem
\ref{bigprodthm} implies these differential operators, between
weighted tensor or spinor bundles, will also factor in a way
generalising Theorem \ref{prodthm}.

The second direction of
generalisation of Theorem \ref{prodthm} concerns the order.  Before
we discuss this we need some definitions.  In the setting of
(pseudo-)Riemannian geometry, we say that a differential operator
between density or tensor bundles is a {\em natural differential
operator} if it can be written as a universal polynomial in
(Levi-Civita connection) covariant derivatives with coefficients
depending polynomially on the metric, its inverse, the curvature
tensor and its covariant derivatives.  The coefficients of natural
operators are called {\em natural tensors}. In the case that they are
scalar they are often also called {\em Riemannian invariants}.  We say
that a differential operator is a {\em conformally invariant
differential operator} if it is a natural operator in this way and is
well-defined on conformal structures (i.e.\ is independent of a choice
of conformal scale). Sometimes for emphasis we shall say such
operators are natural conformally invariant differential operators.
As mentioned above, in even dimensions natural conformally invariant
operators with leading term a power of the Laplacian only exist to
order at most the dimension. We show that on conformally Einstein even
dimensional manifolds the GJMS operators, and their generalisations
(as in Theorem \ref{powerslap}), are extended to all orders by the
operators $P_{k}$ defined in Section \ref{lapsect}.  A qualification
is required: The extending operators $P_{k}$ are canonical and
well-defined (i.e.\ they depend only on the conformal class) on any
conformally Einstein manifold, but at high order on even manifolds
they are not natural conformally invariant differential operators
according to the definition above.

Explaining this last subtle point brings us to the plan and strategy
of the paper. Each Einstein metric corresponds (in a one-one way) to a
suitably generic parallel standard tractor field \cite{Gau,GoNur} that
we denote $\Pa^g$ (or simply $\Pa$ if $g$ is understood). Thus on a
conformally Einstein manifold one has (at least one) parallel tractor
field. We can paraphrase the main construction as follows (and see the
remark after the proof of Theorem \ref{bigprodthm} for a
qualification). The parallel tractor $\Pa^A$ combines with a basic
tool of conformal geometry, the conformally invariant second order
differential tractor D-operator $D_A$, to give a fundamental Laplacian
$\Pa^AD_A$ corresponding to each Einstein metric. Taking this
Laplacian to the power $k$, for each positive integer $k$, yields
differential operators $P_k$ of the form
$\Delta^k+\textit{lower~order~terms}$. Each $P_k$ acts
on tractor fields of weight $k-n/2$ and of any rank (including
densities as the low rank extreme). Thus, essentially by construction,
the operators $P_{k}$ factor along the lines of Theorem
\ref{prodthm}. This is Theorem \ref{bigprodthm}.  It would seem from
this construction that each $P_{k}$ should depend on the particular
Einstein metric $g$ through the corresponding parallel tractor
$\Pa^g$. Remarkably it turns out that if $g_1$ and $g_2$ are
conformally related Einstein metrics then they determine the same
operator $P_{k}$. See Theorem \ref{g1g2}. Thus $P_{k}$ is a canonical
object on conformally Einstein structures in that it depends only on
the conformal class of conformally Einstein metrics.  This surprising
result is partly explained by Theorem \ref{natthm}.  The GJMS
operators $\Box_k$ and their generalisations $\Box_k^0$, as in Theorem
\ref{powerslap}, are natural conformally invariant (density-valued)
operators and therefore if $g_1$ and $g_2$ are conformally related
then $(\Box_k^0)^{g_1}=(\Box_k^0)^{g_2}$. (Naturality for operators
between tractor bundles is defined in section \ref{lapsect}.)
 Theorem \ref{powerslap} shows that when the operator
$\Box_k^0$ exists (i.e.\ the dimension $n$ is odd or, if $n$ even, is
sufficiently large relative to $k$) then on conformally Einstein
manifolds $\Box_k^0$ agrees with $P_{k}$.  Hence the claim that the
$P_{k}$ extend the operators $\Box_k^0$.  In the cases of even
dimensions and high $k$, where the operator $\Box_k^0$ is not defined,
one might conjecture that the operator $P_{k}$ is given by a natural
formula which, although not conformally invariant in general, is
conformally invariant on conformally Einstein structures.  In fact
this is not the case and the situation is more subtle.  In section
\ref{farside} we treat in detail the operator $P_3$ on densities.  We
see there that the existence of $P_3$ in dimension 4 is linked to a
delicate identity on conformally Einstein manifolds which relates the
Bach tensor to the Cotton tensor through the gradient of an Einstein
scale.  There is no obvious way to eliminate the dependence on the
scale albeit that the construction is independent of which conformally
related Einstein scale we choose. (Note that, as Theorem
\ref{bigprodthm} shows, each $P_k$ can be expressed by a natural
formula in terms of an Einstein metric. However this formula is not
generally conformally invariant.)  In summary the main results of this
article are Theorem \ref{Qthm}, Theorem \ref{g1g2}, Theorem
\ref{natthm} and Theorem \ref{bigprodthm}.  Note that in view of
Theorem \ref{natthm}, Theorem \ref{prodthm} above is a special case of
Theorem \ref{bigprodthm}.

In section \ref{background} we review the basic tractor approach to
conformal geometry before introducing some results relating the
conformal Killing field equation to tractor equations. See in
particular Lemma \ref{cKilltractorlem} and Proposition
\ref{cKilltractor}. The latter of these plays a crucial role in the
subsequent study of conformally related Einstein metrics. See Theorem
\ref{2ein2ckill} and Proposition \ref{kdotCetc}.

Acknowledgements. The main ideas for this article were conceived at
the American Institute of Mathematics 2003 workshop ``Conformal
Structure in Geometry, Analysis, and Physics'' and were inspired by
certain presentations and discussions at that meeting. In particular
for the special case of conformally flat Einstein metrics Tom Branson
described a ``spectrum generating'' route to formulae of the form
\nn{prodformula} and \nn{Qformula}.  Robin Graham presented an
alternative means of getting to the same results, using stereographic
projection, and also sketched an argument, using the ambient metric
and naturality arguments, that the same formulae must hold for any
Einstein metric. Graham also indicated that in even dimensions and on
conformally Einstein manifolds there might be an extension of the GJMS
operators to operators of all even orders.  The proofs given here (and
the treatment of generalisations) are based on a distinct idea that
does not require treating conformally flat metrics as a first case.
Proposition \ref{cKilltractor} and some of the other results in
Sections \ref{ceinm} and \ref{multi} were proved at the Spring 2001
session of the Mathematical Sciences Research Institute,
Berkeley. This support by the MSRI is also greatly
appreciated. Conversations with Andreas \v Cap, Michael Eastwood and
Claude LeBrun have also been helpful. Finally I would like to thank
the referee for making several valuable suggestions.

\section{Einstein metrics and conformal geometry}\label{background}

 \subsection{Conformal geometry and tractor calculus}\label{tractorsect}

We first sketch here notation and background for conformal structures.
Further details may be found in \cite{CapGoamb,GoPetLap}.  Let $M$ be
a smooth manifold of dimension $n\geq 3$. Recall that a {\em conformal
structure\/} of signature $(p,q)$ on $M$ is a smooth ray subbundle
$\cq\subset S^2T^*M$ whose fibre over $x$ consists of conformally
related signature-$(p,q)$ metrics at the point $x$. Sections of $\cq$
are metrics $g$ on $M$. So we may equivalently view the conformal
structure as the equivalence class $[g]$ of these conformally related
metrics.  The principal bundle $\pi:\cq\to M$ has structure group
$\Bbb R_+$, and so each representation ${\Bbb R}_+ \ni x\mapsto
x^{-w/2}\in {\rm End}(\Bbb R)$ induces a natural line bundle on $
(M,[g])$ that we term the conformal density bundle $E[w]$. We shall
write $ \ce[w]$ for the space of sections of this bundle. Here and
throughout, sections, tensors, and functions are always smooth.  When
no confusion is likely to arise, we will use the same notation for a
bundle and its section space.

We write $\bg$ for the {\em conformal metric}, that is the
tautological section of $S^2T^*M\otimes E[2]$ determined by the
conformal structure. This will be used to identify $TM$ with
$T^*M[2]$.
For many calculations we will use abstract indices in an obvious way.  
Given a choice of
metric $ g$ from the conformal class, we write $ \nabla$ for the
corresponding Levi-Civita connection. With these conventions the
Laplacian $ \Delta$ is given by $\Delta=\bg^{ab}\nd_a\nd_b=
\nd^b\nd_b\,$.  
Note $E[w]$ is trivialised by a choice of metric $g$ from the
conformal class, and we write $\nd$ for the connection corresponding
to this trivialisation.
It follows immediately that (the coupled) $ \nd_a$  preserves the 
conformal metric.  

Since the Levi-Civita connection is torsion-free, its curvature
$R_{ab}{}^c{}_d$ (the Riemannian
curvature) is given by $ [\nd_a,\nd_b]v^c=R_{ab}{}^c{}_dv^d $
($[\cdot,\cdot]$ indicates the commutator bracket).  The Riemannian
curvature can be decomposed into the totally trace-free Weyl curvature
$C_{abcd}$ and a remaining part described by the symmetric {\em
Schouten tensor} $\Rho_{ab}$, according to $
R_{abcd}=C_{abcd}+2\bg_{c[a}\Rho_{b]d}+2\bg_{d[b}\Rho_{a]c}, $ where
$[\cdots]$ indicates antisymmetrisation over the enclosed indices.
The Schouten tensor is a trace modification of the Ricci tensor
$\Ric_{ab}$ and vice versa: $\Ric_{ab}=(n-2)\Rho_{ab}+\J\bg_{ab}$,
where we write $ \J$ for the trace $ \V_a{}^{a}$ of $ \V$.  The {\em
Cotton tensor} is defined by
$$
A_{abc}:=2\nabla_{[b}\Rho_{c]a} .
$$
Under a {\em conformal transformation} we replace a choice of metric $
g$ by the metric $ \hat{g}=e^{2\om} g$, where $\omega$ is a smooth
function. Explicit formulae for the corresponding transformation of
the Levi-Civita connection and its curvatures are given in e.g.\ 
\cite{BEGo,GoPetLap}. We recall that, in particular, the Weyl curvature is
conformally invariant $\widehat{C}_{abcd}=C_{abcd}$.

We next define the standard tractor bundle over $(M,[g])$.
It is a vector bundle of rank $n+2$ defined, for each $g\in[g]$,
by  $[\ce^A]_g=\ce[1]\oplus\ce_a[1]\oplus\ce[-1]$. 
If $\wh g=e^{2\Up}g$, we identify  
 $(\alpha,\mu_a,\tau)\in[\ce^A]_g$ with
$(\wh\alpha,\wh\mu_a,\wh\tau)\in[\ce^A]_{\wh g}$
by the transformation
\begin{equation}\label{transf-tractor}
 \begin{pmatrix}
 \wh\alpha\\ \wh\mu_a\\ \wh\tau
 \end{pmatrix}=
 \begin{pmatrix}
 1 & 0& 0\\
 \Up_a&\delta_a{}^b&0\\
- \tfrac{1}{2}\Up_c\Up^c &-\Up^b& 1
 \end{pmatrix} 
 \begin{pmatrix}
 \alpha\\ \mu_b\\ \tau
 \end{pmatrix} ,
\end{equation}
where $\Up_a:=\nd_a \Up$.
It is straightforward to verify that these identifications are
consistent upon changing to a third metric from the conformal class,
and so taking the quotient by this equivalence relation defines the
{\em standard tractor bundle} $\ce^A$ over the conformal manifold.
(Alternatively the standard tractor bundle may be constructed as a
canonical quotient of a certain 2-jet bundle or as an associated
bundle to the normal conformal Cartan bundle \cite{luminy}.) On a conformal structure of signature $(p,q)$, the
bundle $\ce^A$ admits an invariant metric $ h_{AB}$ of signature
$(p+1,q+1)$ and an invariant connection, which we shall also denote by
$\nabla_a$, preserving $h_{AB}$.  In a conformal scale $g$, these are
given by
\begin{equation}\label{basictrf}
 h_{AB}=\begin{pmatrix}
 0 & 0& 1\\
 0&\bg_{ab}&0\\
1 & 0 & 0
 \end{pmatrix}
\text{ and }
\nabla_a\begin{pmatrix}
 \alpha\\ \mu_b\\ \tau
 \end{pmatrix}
 =
\begin{pmatrix}
 \nabla_a \alpha-\mu_a \\
 \nabla_a \mu_b+ \bg_{ab} \tau +\Rho_{ab}\alpha \\
 \nabla_a \tau - \Rho_{ab}\mu^b  \end{pmatrix}. 
\end{equation}
It is readily verified that both of these are conformally well-defined,
i.e., independent of the choice of a metric $g\in [g]$.  Note that
$h_{AB}$ defines a section of $\ce_{AB}=\ce_A\otimes\ce_B$, where
$\ce_A$ is the dual bundle of $\ce^A$. Hence we may use $h_{AB}$ and
its inverse $h^{AB}$ to raise or lower indices of $\ce_A$, $\ce^A$ and
their tensor products.

In computations, it is often useful to introduce 
the `projectors' from $\ce^A$ to
the components $\ce[1]$, $\ce_a[1]$ and $\ce[-1]$ which are determined
by a choice of scale.
They are respectively denoted by $X_A\in\ce_A[1]$, 
$Z_{Aa}\in\ce_{Aa}[1]$ and $Y_A\in\ce_A[-1]$, where
 $\ce_{Aa}[w]=\ce_A\otimes\ce_a\otimes\ce[w]$, etc.
 Using the metrics $h_{AB}$ and $\bg_{ab}$ to raise indices,
we define $X^A, Z^{Aa}, Y^A$. Then we
immediately see that 
$$
Y_AX^A=1,\ \ Z_{Ab}Z^A{}_c=\bg_{bc} ,
$$
and that all other quadratic combinations that contract the tractor
index vanish. 
In \eqref{transf-tractor} note that  
$\wh{\alpha}=\alpha$ and hence $X^A$ is conformally invariant.

Given a choice of conformal scale, the {\em tractor-$D$ operator} 
$$
D_A\colon\ce_{B \cdots E}[w]\to\ce_{AB\cdots E}[w-1]
$$
is defined by 
\begin{equation}\label{Dform}
D_A V:=(n+2w-2)w Y_A V+ (n+2w-2)Z_{Aa}\nabla^a V -X_A\Box V, 
\end{equation} 
 where $\Box V :=\Delta V+w \J V$.  This also turns out to be
 conformally invariant as can be checked directly using the formulae
 above (or alternatively there are conformally invariant constructions
 of $D$, see e.g.\ \cite{Gosrni}).

The curvature $ \Omega$ of the tractor connection 
is defined by 
$$
[\nd_a,\nd_b] V^C= \Omega_{ab}{}^C{}_EV^E 
$$
for $ V^C\in\ce^C$.  Using
\eqref{basictrf} and the formulae for the Riemannian curvature yields
\begin{equation}\label{tractcurv}
\Omega_{abCE}= Z_C{}^cZ_E{}^e C_{abce}-2X_{[C}Z_{E]}{}^e A_{eab}
\end{equation}

We will also need a conformally invariant curvature quantity defined as
follows (cf.\ \cite{Gosrni,Goadv})
\begin{equation}\label{Wdef}
W_{BC}{}^E{}_F:=
\frac{3}{n-2}D^AX_{[A} \Omega_{BC]}{}^E{}_F ,
\end{equation}
where $\Omega_{BC}{}^E{}_F:= Z_A{}^aZ_B{}^b \Om_{bc}{}^E{}_F$.
In a choice of conformal scale, 
 $W_{ABCE}$ is given by
\begin{equation}\label{Wform}
\begin{array}{l}
(n-4)\left( Z_A{}^aZ_B{}^bZ_C{}^cZ_E{}^e C_{abce}
-2 Z_A{}^aZ_B{}^bX_{[C}Z_{E]}{}^e A_{eab}\right. \\ 
\left.-2 X_{[A}Z_{B]}{}^b Z_C{}^cZ_E{}^e A_{bce} \right)
+ 4 X_{[A}Z_{B]}{}^b X_{[C} Z_{E]}{}^e B_{eb},
\end{array}
\end{equation}
where 
$$
B_{ab}:=\nabla^c
A_{acb}+\Rho^{dc}C_{dacb}.
$$ is known as the {\em Bach tensor}. From the formula \nn{Wform} it
is clear that $W_{ABCD}$ has Weyl tensor type symmetries.  Another way
to find the tractor  $W_{ABCD}$ is through commuting tractor D-operators: for 
$V^A\in \ce^A$ we have
\begin{equation}\label{Dcomm} 
[D_A,D_B]V^K=(n+2w-2)W_{AB}{}^{K}{}_{L}V^L+6X_{[A}
\Omega_{BP]}{}^{K}{}_{L}D^PV^L ~.
\end{equation}
 The tractor field $W_{ABCD}$ has an important
relationship to the ambient metric of Fefferman and Graham.  For a
conformal manifold of signature $(p,q)$ the ambient manifold
\cite{FGast} is a signature $(p+1,q+1)$ pseudo-Riemannian manifold
with $\cq$ as an embedded submanifold. Suitably homogeneous tensor
fields on the ambient manifold, upon restriction to $\cq$, determine
tractor fields on the underlying conformal manifold \cite{CapGoamb}.
In particular, in dimensions other than 4, $W_{ABCD}$ is the tractor
field equivalent to $(n-4){\aR}|_\cq$ where $\aR$ is the curvature of
the Fefferman-Graham ambient metric.

Many natural conformal equations have elegant and useful
interpretation in terms of tractors. A case we will need below is the
conformal Killing equation. Recall that a tangent field $k$ is a
conformal Killing field if, for any metric $g$ in
the conformal class, there is a function $\lambda$ so that  
$\cL_k g=\lambda g$. In terms of the Levi-Civita connection for $g$ this equation is given by
$$
\nd_{(a}k_{b)_0}=0 ,
$$
where $(\cdots)_0$ indicates the trace-free symmetric part.
To find the equivalent tractor expression we 
first note that, as observed in \cite{Esrni}, 
the operator $M_{A}{}^a:\ce_a[w]\to \ce_A[w-1]$ given, in a conformal
scale by
$$
u_a\mapsto (n+w-2)Z_A{}^a u_a -X_A\nd^a u_a,
$$
is conformally invariant. We obtain the following: 
\begin{lemma}\label{cKilltractorlem}
Solutions $k^a$ of the conformal Killing equation
$$
\nd_{(a}k_{b)_0}=0
$$
are in 1-1 correspondence with fields $K_B\in\ce_B[1]$ 
such that  
\begin{equation}\label{cksttr}
D_{(A}K_{B)}=0 
\end{equation} 
The correspondence is given by 
$$ k^a\mapsto \frac{1}{n}M_{Ba}k^a \quad \mbox{with inverse} \quad K_{B}
\mapsto Z^{Ba}K_{B}.
$$  
\end{lemma}
\noindent{\bf Proof:} Suppose that $K_B\in\ce_B[1] $ satisfies
$D_{(A}K_{B)}=0$. Then $0=X^AX^B D_{A}K_{B}= nX^B K_B$. So $K_B$ takes 
the form $Z_B{}^b k_b -X_B \rho$ for some tangent field $k^b$ and function $\rho$.  
But, using this, a short calculation shows that $D^BK_B=0$ implies that $\rho =\nd^a k_a/n$, that is 
\begin{equation}\label{KKform}
K_B=  Z_B{}^b k_b - \frac{1}{n}X_B\nd^a k_a =  \frac{1}{n}M_{Ba}k^a.
\end{equation}
Using this, and \nn{basictrf}, \nn{Dform} once again, we calculate that 
\begin{equation}\label{big}
\begin{aligned}
&D_{(A} K_{B)}=n^2 \nd_{(a}k_{b)_0}Z_A{}^aZ_B{}^b\\
&+[(2-n)(\nd_b\nd^ck_c+nP_{bc}k^c)-n\Delta k_b-nJ k_b]X_{(A}Z_{B)}{}^b\\
&-[\Delta \nd^ck_c+nk^c\nd^bP_{bc} +nP_{bc}\nd^bk^c +J \nd_ck^c ]X_AX_B .
\end{aligned}
\end{equation}
 Thus $D_{(A}K_{B)}=0$ implies the vanishing of
$\nd_{(a}k_{b)_0}$. 

Now let us assume that $\nd_{(a}k_{b)_0}=0$. It remains to show that
the last display vanishes. We have that 
\begin{equation}\label{ndk}
\nd_a k_b=\mu_{ab}+\nu
\bg_{ab}
\end{equation}
 where $\mu_{ab}=-\mu_{ba}$. Differentiating and using the
Bianchi identity reveals that 
$$
\nd_{a}\mu_{bc}=\nd_{c}\mu_{ba}-\nd_b\mu_{ca} .
$$
Substituting, on the right-hand side, for $\mu$ using \nn{ndk}, and
using the formula for the Riemannian curvature in terms of the Weyl
curvature and the Schouten tensor, we obtain
\begin{equation}\label{ndmu}
\nd_a \mu_{bc}=\bg_{ab}\rho_c-\bg_{ac}\rho_b-P_{ab}k_c 
+P_{ac}k_b +C^{\phantom{bc}}_{bc}{}^d{}^{\phantom{d}}_ak_d ,
\end{equation}
where $\rho_c=-P_c{}^b k_b -\nd_c \nu$. We use this result to replace 
$\nd_a\mu_{cb}$ in the expansion
$$
\Delta k_b=\bg^{ac}\nd_a\nd_c k_b=\bg^{ac}\nd_a\mu_{cb}+\nd_b\nu , 
$$
to obtain
$$
\Delta k_b =(n-2)\rho_b -J k_b. 
$$
On the other hand, from the definitions of $\nu$ and $\rho_c$ we have 
$$
\nd_b \nd_c k^c =-n \rho_b -n P_{bc}k^c .
$$
Together these show that in the display \nn{big} the 
coefficient of $X_{(A}Z_{B)}{}^b$ vanishes. 

Next, contracting \nn{ndmu} with $\bg^{ab}$ and then taking a
divergence establishes that $\nd^b\rho_b=J\nu$. Using this it follows
easily that $\Delta \nd_b k^b= -2n J \nu - n k^b(\nd_b J)$ (after
using a contracted Bianchi identity and \nn{ndk}).  Substituting this
into \nn{big} we see that the coefficient of
$X_AX_B$ also vanishes.  \quad $\Box$\\
\noindent{\bf Remark:} The main result of the Lemma above can be seen
to hold without performing the explicit calculations described. We
sketch the idea. The differential operator $M_{Ba}$ given above is a
differential splitting operator. It is easily verified explicitly that
there is also a differential splitting operator
$S^{ab}_{AB}\ce_{(ab)_0}[2]\to \ce_{(AB)}$.  From the classification
of conformally invariant operators on the conformal sphere (see e.g.\
\cite{ESlo}) and the composition series for $\ce_{(AB)}$ it follows
easily that on the sphere $D_{(A}M_{B)}{}^a k_a -
S^{ab}_{AB}\nd_{(a}k_{b)_0}=0$. In fact this also holds on conformally
curved structures. Clearly the difference must involve curvature. From
the conformal invariance it is not difficult to argue that the
identity could only fail to hold if there is some partial contraction
involving at most one use of $\nd$ and otherwise involving only (and
non-trivially) one power of the Weyl curvature, and also $k$ linearly,
and such that this partial contraction could appear in the composition
series for $\ce_{(AB)}$. From elementary invariant theory and weight
considerations it is easily verified that there is no such partial
contraction.

\medskip

\noindent The result in the Lemma should be of some independent interest since,
via the machinery of \cite{CapGoamb,GoPetLap}, it yields an ambient
metric interpretation of infinitesimal conformal symmetries. Our
immediate interest is that it leads to a simple and suggestive proof
of the following result.
\begin{proposition}\label{cKilltractor}
Solutions $k^a$ of the conformal Killing equation
$$
\nd_{(a}k_{b)_0}=0
$$
are in 1-1 correspondence with solutions $\K_{DE}\in \ce_{[DE]}$ of the 
equation
\begin{equation}\label{cktreq}
\nd_a \K_{DE}-\K_{AB}X^AZ^{Ba}\Omega_{abDE}=0.
\end{equation}
The correspondence is given by 
$$ k^a\mapsto \frac{1}{n^2} D_{[A}M_{B]a}k^a \quad \mbox{with inverse} 
\quad \K_{AB}
\mapsto X^AZ^{Ba} \K_{AB}.
$$  

If $\K_{AB}$ is a parallel adjoint tractor then $k^a:=X^AZ^{Ba}
\K_{AB}$ is a conformal Killing vector field and hence, from \nn{cktreq}, satisfies $k^a\Omega_{abDE}=0 $.
\end{proposition}
\noindent{\bf Proof:} Suppose that $\K_{DE}$ solves \nn{cktreq} or
alternatively that $\K_{DE}$ is parallel. Since $X^D \Omega_{abDE}=0$
it follows that in either case
$$
X^BZ^C{}_c \nd_a\K_{BC}=0 .
$$ Expanding this out using \nn{basictrf} etcetera and writing
$k_c:= X^BZ^C{}_c \K_{BC}$ we obtain
$$
\nd_a k_c=\mu_{ac} + f \bg_{ac}
$$
for some function $f$ and skew 2-form 
$\mu$ of conformal weight $2$.

Now for the opposite implication, suppose that $k^a$ is a conformal Killing vector field. 
From the Lemma $K_A:=\frac{1}{n}M_A{}^ak_a$ satisfies 
$$
D_A K_B=-D_B K_A = n\K_{AB},
$$
where we have written $\K_{AB}$ as a shorthand for $\frac{1}{n^2}
D_{[A}M_{B]a}k^a$.
This implies 
$$
2nD_A\K_{BC}=2D_A D_B K_C =[D_A,D_B]K_C+[D_C,D_A]K_B + [D_C,D_B]K_A.
$$ Now note that, since $\K_{BC}$ has conformal weight 0, it follows
immediately from \nn{Dform} that $X^AD_A\K_{BC}=0 $. Observe also that  $Z^A{}_a
2nD_A\K_{BC}$ recovers $ 2n(n-2)\nd_a\K_{BC}$. 
Next from
\nn{Dcomm} we have 
\begin{eqnarray*} 
[D_A,D_B]K_C &=&nW_{ABCE}K^E+6X_{[A}
\Omega_{BP]CE}D^P K^E\\
&=&nW_{ABCE}K^E+6nX_{[A}
\Omega_{BP]CE} \K^{PE} .
\end{eqnarray*}
Now let us
consider first dimensions $n\neq 4$. In this case we have $X_{[A}
\Omega_{BP]CE} = \frac{1}{n-4} X_{[A} W_{BP]CE} $, from \nn{Wform}, 
and so substituting this into the last but one display, we have 
$$
\begin{aligned}
2nD_A\K_{BC}=& 2n W_{CBAQ}K^Q \\
&+\frac{6n}{n-4}(X_{[C}W_{BS]AQ}+X_{[A}W_{BS]CQ}+ X_{[C}W_{AS]BQ})\K^{SQ} ,
\end{aligned}
$$ where we have used that $W_{[BCA]Q}=0$ to simplify.  Expanding the
skew symmetrisations, in the display, we find that 4 terms cancel due
to the Weyl-tensor like symmetries of the W-tractor (in particular
that these imply that e.g.\ $W_{ASBQ}\K^{SQ}$ is skew). Contracting
both sides with $Z^A{}_a$ and dividing by $2n$ brings us to
$$
(n-2) \nd_a\K_{BC} = Z^A{}_a (W_{CBAQ}K^Q + \frac{2}{n-4}X_S W_{CBAQ}\K^{SQ}).
$$ These terms combine as $X_S\K^{SQ} = K^Q$.  But from \nn{Wform},
\nn{KKform} and the symmetries of $W$ we have $
W_{CBAQ}K^QZ^A{}_a=(n-4)k^q\Omega_{qaBC}$ and so the right-hand side
is exactly $(n-2)k^q\Omega_{qaBC}$ as required; that is we obtain
\nn{cktreq}.  For the calculation in dimension 4 the manoeuvre of
replacing $X_{[A} \Omega_{BP]CE} = \frac{1}{n-4} X_{[A} W_{BP]CE} $ is
not possible. In dimension 4 contraction with $Z^A{}_a$ annihilates
$W_{CBAQ}$ and we come to
$$
\nd_a\K_{BC} = 3Z^A{}_a(X_{[C}\Omega_{BS]AQ}+X_{[A}\Omega_{BS]CQ}+ X_{[C}\Omega_{AS]BQ})\K^{SQ} .
$$ 
Expanding this using \nn{tractcurv}, a
straightforward calculation and use of the symmetries of the Weyl tensor
yields \nn{cktreq}.
 \quad $\Box$\\
\noindent{\bf Remark:} The idea of the proof above is clarified
by considering the conformally flat case.  On conformally flat
structures we have $[D_A,D_B]=0$ when acting on any tractor so $ D_A
D_B K_C  $ is symmetric on the index pair
``{\scriptsize$AB$}''. On the other hand from the Lemma it is skew on
the pair index pair ``{\scriptsize$BC$}''. Together these imply 
$$ D_A
D_B K_C =0 \quad \Rightarrow \quad \nd_a D_B K_C=n\nd_a \K_{BC} =0 . 
$$ 

 Finally note that since the term $\K_{AB}X^AZ^{Ba}\Omega_{abDE}$, on
the left hand side of \nn{cktreq}, is linear in $\K_{AB}$, it follows
that $X^{[A}Z^{B]a}\Omega_{abDE}$ may be viewed as a contorsion that
modifies the tractor connection to a new conformally invariant
connection $\tilde{\nd}$ on the bundle $\ce_{[AB]}$.
Then sections of $\ce_{[AB]}$ which are parallel for $\tilde{\nd}$ are in 1-1 correspondence with  conformal
Killing vector fields.  \quad \endrk

\subsection{Conformally Einstein manifolds}\label{ceinm}

Recall that  a conformal structure $[g]$ is said to be conformally
Einstein if there is a metric $\widehat{g}$ in the conformal class
(i.e. $\widehat{g}\in [g]$) such that the Schouten tensor for
$\widehat{g}$ is pure trace.  Let us use the term {\em Einstein
tractor} for a parallel standard tractor $\Pa$ with the property that
$X_A\Pa^A$ is nowhere vanishing. We have the
following result.
\begin{proposition}\cite{Gau,GoNur}
\label{cein}
On a conformal manifold $(M,[g])$ there is a 1-1 correspondence
between conformal scales $\si\in \ce[1]$, such that
$g^\si=\si^{-2}\bg$ is Einstein, and Einstein tractors. The mapping
from Einstein scales to parallel tractors is given by 
$\si\mapsto \frac{1}{n}D_A \si$ while the inverse is $\Pa^A \mapsto X_A\Pa^A$.
\end{proposition}

It follows that if $g$ is an Einstein metric and $\Pa$ the
corresponding tractor then
$\Omega_{bc}{}^D{}_E\Pa^E=[\nd_b,\nd_c]\Pa^D=0$. Also since $\Pa$ is
parallel and of weight 0 then, viewing it as a multiplication
operator, we have $[D,\Pa]=0$.  {}From \nn{Wdef} we have $
W_{A_1A_2}{}^D{}_E = \frac{3}{n-2}D^{A_0}X_{A_0}
Z_{A_1}{}^bZ_{A_2}{}^c\Omega_{bc}{}^D{}_E$.  Thus $W_{BCDE}\Pa^E =0 $.
On the other hand from \nn{Wform} it follows that $W_{BCDE}$ has Weyl
tensor type symmetries. Thus, in summary:
\begin{equation}\label{PaW}
\Omega_{bc}{}^D{}_E\Pa^E=0, \quad W_{BCDE}\Pa^E =0, \quad {\rm and}
\quad \Pa^B W_{BCDE}=0 .
\end{equation}

\subsection{Conformally related Einstein metrics}\label{multi}

There is an intimate relationship between pairs of conformally related
Einstein metrics and conformal Killing fields. A part of this story goes
as follows.
\begin{theorem}\label{2ein2ckill}
  If $g_1 = (\si_1)^{-2}\bg$ is an Einstein metric, then the
  conformally related metric $g_2= (\si_2)^{-2}\bg $ is an Einstein
  metric if and only if the vector field
\begin{equation}\label{kdef}
 k^a:=\si_1\nd^a\si_2 -\si_2\nd^a \si_1
\end{equation}
is a conformal Killing vector field. 
\end{theorem}
\noindent{\bf Proof:} $\Rightarrow:$
 Since $g_1$ and $g_2$ are Einstein it follows from Proposition
\ref{cein} that 
$$
\Pa_1^A := \frac{1}{n} D^A \si_1 \quad \mbox{ and } \Pa_2^A:= \frac{1}{n}D^A\si_2 
$$
are both parallel for the tractor connection. Thus 
\begin{equation}\label{kpar}
\K^{AB}:= \Pa_1^A\Pa_2^B- \Pa_1^B\Pa_2^A
\end{equation} 
is a parallel adjoint tractor. But, from the formula \nn{Dform} for
$D$, this has primary part (i.e.\ $k^b:=X_AZ_{B}{}^b\K^{AB}$) given by
$$
k_b=\si_1\nd_b\si_2 -\si_2\nd_b \si_1
$$ where $\nd_b$ is the Levi-Civita connection for any metric $g$ in
the conformal class $[g_1]=[g_2]$. By Proposition \ref{cKilltractor} $k^a$ is
a conformal Killing field.

$\Leftarrow :$ From expression (4), in section 2.3 of \cite{BEGo}, we
have that $g_2$ is Einstein if and only if
$\nd_{(a}^{g_1}\nd^{g_1}_{b)_0}\si_2 + P^{g_1}_{(ab)_0}\si_2=0$. Now recall
that $\nd^{g_1}_a\si_1=0 $. Using this and that $\si_1$ is
non-vanishing, we see that $\nd_{(a}k_{b)_0}=0$ implies
$\nd_{(a}^{g_1}\nd^{g_1}_{b)_0}\si_2=0$.  On the other hand since
$g_1$ is Einstein $P^{g_1}_{(ab)_0}=0$. \quad $\Box$

\medskip

\noindent The Theorem above is classical and probably due to Brinkmann
\cite{Brinkman}.  Our treatment via tractors is helpful for arguments
which follow. We will also need the following.
\begin{proposition}\label{kdotCetc} If the conformally related metrics 
$g_1 = (\si_1)^{-2}\bg$  and 
 $g_2= (\si_2)^{-2}\bg $ are both Einstein  then 
$$ k^a\Om_{ab}{}^C{}_D =0 \quad {\rm and} \quad (\nd^a
\si_1)(\nd^b\si_2)\Om_{ab}{}^C{}_D=0,
$$
where $k^a$ is the the conformal Killing vector field
$k^a$, as given in \nn{kdef}.
\end{proposition}
\noindent{\bf Proof:} Let $\K$ be defined as in \nn{kpar}. Since $\K$
is parallel it follows immediately from Proposition \ref{cKilltractor}
that $k^a$ satisfies
$$
k^a\Omega_{abCD}=0.
$$

Next, using the shorthand $\mu^a_i:=\nd^a \si_i$, $i=1,2$, we have
$k^a=\si_1\mu_2^a-\si_2\mu_1^a$ and 
$$
k^a \Omega_{abCD}=0  \Rightarrow \si_1\mu_2^a \Omega_{abCD} =
\si_2 \mu_1^a\Omega_{abCD},
$$
and so 
$$
 \si_1\mu_1^b \mu_2^a \Omega_{abCD}  =\si_2\mu_1^b\mu_1^a \Omega_{abCD}=0,
$$
since $\Omega_{abCD}$ is skew on the index pair ``\mbox{\scriptsize{$ab$}}''.
Thus 
$$
\mu_1^a \mu_2^b \Omega_{abCD}=0,
$$
as claimed.
\quad $\Box$\\
\noindent{\bf Remark:} Note that via \nn{Dform} we can state the
results of the Proposition as follows: for Einstein
tractors $\Pa_1$ and $\Pa_2$ we have 
$\Pa_1^A\Pa_2^B Z_A{}^aZ_B{}^b\Omega_{ab}{}^C{}_D=0 $ and
$\Pa_1^A\Pa_2^BX_{[A}Z_{B]}{}^b \Omega_{bc}{}^D{}_E=0 $.  \quad \endrk

\section{Conformal operators of Laplace type}\label{lapsect}

Let us write $\ce^\bullet[w]$ to indicate an any fixed tractor bundle or
density bundle of weight $w$.\\
\noindent{\bf Definition:} On a conformally Einstein manifold with
$g=\si^{-2}\bg$ an Einstein metric we define, for each positive
integer $k$, a differential operator
$$
P^{g}_{k}:\ce^\bullet[k-n/2]\to \ce^\bullet[-k-n/2]
$$
by 
\begin{equation}\label{Pdef}
P^{g}_{k} f=(-1)^{k-1}\si^{1-k} \Pa^{A_1}\cdots \Pa^{A_{k-1}} \Box D_{A_1}\cdots D_{A_{k-1}} f ,
\end{equation}
where $\Pa_A=\frac{1}{n}D_A\si$. It follows easily from \nn{Dform}
 that this has the form $P^{g}_{k}=\Delta^k+${\em lower order
 terms}.

Since the tractor D-operator is conformally invariant on weighted
tractor bundles and $\Box$ is conformally invariant on tractor bundles
of weight $(1-n/2)$ it follows that the composition $\Box
D_{A_1}\cdots D_{A_{k-1}} $ is conformally invariant on
$\ce^\bullet[k-n/2]$. Through the tractor field
$\Pa_A=\frac{1}{n}D_A\si$ it would appear that the operator
$P^{g}_{k} $ depends on the choice of a particular Einstein metric
$g =\si^{-2}\bg$ from the conformal class of metrics.  The following
theorem shows that in fact the operator $P^{g}_{k}$ is a canonical
object depending only on the conformal class of conformally Einstein
metrics.
 \begin{theorem}\label{g1g2} If $g_1$ and $g_2$ are two conformally
  related Einstein metrics then  $P^{(g_1)}_{k}= P^{(g_2)}_{k}$
\end{theorem} 
\noindent{\bf Proof:} First note that if $f\in \ce^\bullet[k-n/2]$
then $D_{A_1}\cdots D_{A_{k-1}} f $ has weight $1-n/2$ and so,
from the formula for the tractor operator \nn{Dform}, we have
\begin{equation}\label{Xthere}
D_{A_0}D_{A_1}\cdots D_{A_{k-1}}f=-X_{A_0}\Box D_{A_1}\cdots
D_{A_{k-1}}f.
\end{equation}  
Thus if $\Pa$ is an Einstein tractor 
 corresponding to $g=\si^{-2}\bg$ then we have
\begin{equation}\label{altform}
P^{g}_{k} f=(-1)^{k}\si^{-k} \Pa^{A_0}\Pa^{A_1}\cdots
\Pa^{A_{k-1}} D_{A_0} D_{A_1}\cdots D_{A_{k-1}} f .
\end{equation}

Consider a conformal manifold $(M,\bg)$.
Suppose that for $i=0,1,\cdots ,k$, $\si_i$ are (not necessarily
distinct) Einstein scales for the  conformal structure $\bg$.
 That is
$$
g_i:=\si_i^{-2}\bg \quad \mbox{ is Einstein } \quad i=0,1,\cdots ,k.
$$ 
Write $\Pa_i$ for the parallel tractor corresponding to $\si_i$ (in
the sense of Proposition \ref{cein}) as $i$ ranges over $0,\cdots ,k$.
In view of \nn{altform}, to prove the theorem it obviously suffices to
show that for any $\ell\in \{0,1,\cdots ,k-1 \}$ we have
\begin{equation}\label{equiv}
\begin{aligned}
&\si_0^{-1} \cdots \si_{k-1}^{-1} \Pa_0^{A_0}\cdots
\Pa_{k-1}^{A_{k-1}} D_{A_0} \cdots D_{A_{k-1}} f\\
 & =\tilde{\si}_0^{-1} \cdots \tilde{\si}_{k-1}^{-1}
\tilde{\Pa}_0^{A_0}\cdots
\tilde{\Pa}_{k-1}^{A_{k-1}} D_{A_0} \cdots D_{A_{k-1}} f
\end{aligned}
\end{equation}
where $\tilde{\si}_\ell=\si_k$ and $\tilde{\Pa}_\ell=\Pa_k$ but
otherwise, for $i\in \{0,\cdots,k-1\}\setminus \{\ell \}$ we have
$\tilde{\si}_i=\si_i$ and $\tilde{\Pa}_i=\Pa_i$.  Observe first that
from \nn{Xthere} this is clear in the case $\ell=0$ since $\si_0^{-1}
\Pa_0^{A_0}X_{A_0} = 1= \si_k^{-1}\Pa_k^{A_0}X_{A_0}$.  But this
proves the general case from the separate observation that there is an
equality of the form \nn{equiv} if the tilded quantities on the
right-hand-side are simply a permutation of the untilded quantities on the
left-hand-side.

To establish this last claim it suffices to show that 
$$
\Pa_0^{A_0}\cdots
\Pa_{k-1}^{A_{k-1}} D_{A_0} \cdots D_{A_{\ell-1}}[D_{A_\ell},D_{A_{\ell+1}}]
D_{A_{\ell+2}} \cdots D_{A_{k-1}} f 
$$ vanishes identically for $\ell\in\{ 0,\cdots,k-2 \}$. Now each of
the $\Pa_i$ has weight 0 and is parallel and so commutes with the
tractor-D operators. A straightforward calculation verifies that the
generalisation of \nn{Dcomm} to arbitrary rank tractors $ V\in
\ce_{CE\cdots F}[w]$ may be expressed as follows (cf.\ \cite{GoPetLap})
$$
\begin{array}{lll}
\lefteqn{[D_{A},D_{B}] V_{CE\cdots F}=}&&
\\
&&
(n+2w-2)[W_{ABC}{}^Q V_{QE\cdots F} +2w \Omega_{ABC}{}^Q V_{QE\cdots F}
\\
&&
+4X_{[A}\Omega_{B]}{}^{s}{}_C{}^Q \nd_sV_{QE\cdots F}+\cdots
+ W_{ABF}{}^Q V_{CE\cdots Q}
\\&&
+2w \Omega_{ABF}{}^Q V_{CE\cdots Q}
+4X_{[A}\Omega_{B]}{}^{s}{}_F{}^Q \nd_sV_{CE\cdots Q}] 
,
\end{array}
$$ where $\Omega_{B}{}^{s}{}_C{}^Q:=Z_B{}^b\Omega_{b}{}^s{}_C{}^Q$
and, recall, $ \Omega_{ABC}{}^Q:=Z_A{}^aZ_B{}^b\Omega_{abC}{}^Q$.
Thus the result follows immediately from expression \nn{PaW} and
(the remark following) Proposition \ref{kdotCetc}, which together show that for
any Einstein tractors $\Pa_1$ and $\Pa_2$ we have $\Pa_1^A
W_{ABCD}=0$, $\Pa_1^A\Pa_2^B \Omega_{ABC}{}^Q=0 $ and
$\Pa_1^A\Pa_2^BX_{[A}\Omega_{B]}{}^{s}{}_F{}^Q=0 $.  \quad $\Box$ 

\noindent The Theorem above shows that the operator
$P^{g}_{k}:\ce^\bullet[k-n/2]\to \ce^\bullet[-k-n/2]$ is determined by
the (conformally Einstein) conformal structure; it does not otherwise
depend on any choice of metric from the conformal class.  Thus it is
reasonable to drop reference to a metric in the notation and write
simply $P_{k}$ for this operator.

\noindent{\bf Remarks:} In the case of the theorem where $f$ is simply a 
density there is a simpler conclusion to the proof. 
We may use that $[D,\Pa_i]=0$ to obtain 
$$
\begin{aligned}
&\Pa_0^{A_0}\cdots
\Pa_{k-1}^{A_{k-1}} D_{A_0} \cdots D_{A_{\ell-1}}[D_{A_\ell},D_{A_{\ell+1}}]
D_{A_{\ell+2}} \cdots D_{A_{k-1}} f
\\
&= \Pa_0^{A_0}\cdots
\Pa_{\ell+1}^{A_{\ell+1}}D_{A_0} \cdots D_{A_{\ell-1}}[D_{A_\ell},D_{A_{\ell+1}}](\Pa_{\ell+2}^{A_{\ell+2}}
\cdots
\Pa_{k-1}^{A_{k-1}}
 D_{A_{\ell+2}} \cdots D_{A_{k-1}} f)
\end{aligned}
$$
and then simply use that
$$[D_{A_\ell},D_{A_{\ell+1}}](\Pa_{\ell+2}^{A_{\ell+2}}\cdots
\Pa_{k-1}^{A_{k-1}}
 D_{A_{\ell+2}} \cdots D_{A_{k-1}} f) =0 
$$
as $[D_A,D_B]$ vanishes on densities. 

In another direction, since the theorem holds for tractor valued
densities one might hope that a coupled version of the theorem holds
for densities taking values in some other vector bundle with
connection. One can couple the tractor connection and 
so obtain coupled conformally invariant tractor-D
operators on such bundles and their weighted variants. Note however
that the general proof of the Theorem above does not obviously
generalise to cover such cases since it critically uses the
relationship of the tractor curvature to parallel tractors.  \quad
\endrk

For each choice of conformal scale, the tractor bundles decompose into
direct sums of weighted tensor bundles. So, at each such  scale,
differential operators between (weighted) tractor bundles decompose
into a matrix of differential operators between weighted tensor
bundles. A differential operator between weighted tractor bundles is
said to be a {\em conformally invariant differential operator} if it
is well-defined on conformal structures (i.e.\ is independent of a
choice of conformal scale) and if also, for any metric from the
conformal class, the component operators, in the matrix of operators
between tensor components, are each natural.
Sometimes for emphasis we will term such an operator a  
{\em natural conformally invariant differential operator}). 
 On general
conformal manifolds we have the following result.
\begin{theorem} \label{powerslap}
  For every integer $k\geq 1$ there exist natural conformally
  invariant operators $\Box^0_k:\ce^{\bullet}[k-n/2]\to
  \ce^{\bullet}[-k-n/2]$ having leading term $\Delta^k$ as follows: in
  odd dimensions for all $k\geq 1$; in even dimensions for $1\leq
  k\leq n/2-1$, or if $\ce^{\bullet}[k-n/2]=\ce[k-n/2]$ for $1\leq
  k\leq n/2$. In any dimension if $\ce^{\bullet}[k-n/2]=\ce[k-n/2]$
  then $\Box_k^0$ is a GJMS operator $\Box_k$.

The operators
$\Box_k^{0}$ have tractor formulae as follows:

\begin{equation}\label{powerform}
\begin{array}{lll}
\lefteqn{\displaystyle (-1)^{k-1} X_{A_1}\cdots
X_{A_{k-1}} \Box^0_k U }
&&\vspace{2mm}
\\
&&
\displaystyle
=
\Box D_{A_1}\cdots D_{A_{k-1}}
U
+{\Cal P}^{\bullet,k}_{A_1\cdots A_{k-1}} U
,
\end{array}
\end{equation}
%
where the differential operator ${\Cal P}^{\bullet,k}$ is a partial
contraction polynomial in $X$, $D$, $W$, $h$, and $h^{-1}$. In every
term of this partial contraction at least one of the free indices
$A_1,\cdots ,A_{k-1}$ appears on a $W$.  The coefficients in this
polynomial are rational functions of the dimension.
\end{theorem}
\noindent{\bf Proof:} 
This theorem summarises selected results from
Theorem 4.1 and Theorem 4.14 of \cite{GoPetobstrn} except for the last
two statements. (The notation $\Box_k^0$ is from \cite{GoPetobstrn}.)
The claim concerning rationality is immediate from the
algorithm establishing the proofs of these theorems.

To establish the result concerning the placement of the free indices
let us sketch how get to the formula \nn{powerform}. The operator $\Box_k^0$
arises from the power $\aDe^k$ of the ambient Bochner Laplacian
$\aDe=\aNd^A\aNd_A$ applied to suitably homogeneous ambient
tensors via arguments as in \cite{GoPetLap} (or
\cite{GoPetobstrn}, and the reader should refer to these sources for
further details and the notation). In particular we may follow the procedure established in
\cite{GoPetLap} for the case of densities.\\
Let $\tU$ be the homogeneous (of weight $k-n/2$) ambient tensor
corresponding to the tractor field $U$.  \\
{\bf Step 1:} Observe that
$$\aDe \bD_{A_{k-1}}\cdots \bD_{A_1} \tU =\aDe
(2\aNd_{A_{k-1}}-\X_{A_{k-1}}\aDe)\cdots
(2(k-1)\aNd_{A_1}-\X_{A_1}\aDe)\tU.$$
Expand this out via the
distributive law
without changing the order of any of the operators. \\
{\bf Step 2:} Move all $ \X$'s to the left of any $ \aNd$ or $ \aDe$ via
the identities $ [\aNd_A,\X_B]=\h_{AB}$ and
$ [\aDe , \X_A]= 2\aNd_A$ (which hold to all orders).  \\
{\bf Step 3:} Move all $\aDe$'s to the right of any $\aNd$'s (other
than those implicit in $ \aDe$) via the Ricci identity.  We thus obtain
\begin{equation}\label{form}
        (-1)^{k-1}\X^{k-1}\aDe^k \tU + \sum \h^s \X^x
        (\aNd^{p_1}\aDe^{r_1}\aR) \cdots (\aNd^{p_d}\aDe^{r_d}\aR) \aNd^q
        \aDe^r \tU
        ,
\end{equation}
where $d\geq 1$ in each term of the right-hand part and in each such
term at least one of the $r_d$ occurrences of $\aR$ has one of the
indices $A_1,\cdots , A_{k-1},$ as a free index. This last point is the key 
and follows easily from the following observation: 
If $[\aNd_{A},\aDe]$ vanishes formally then, for $i\in \{ 1,\cdots ,k-1\}$, 
$$
\begin{array}{lll}
\lefteqn{\displaystyle \aDe^{i}\D_{A_{i}}\D_{A_{i+1}}\cdots \D_{A_{k-1}}\tU }
&&\vspace{2mm}
\\
&&
\displaystyle
=\aDe^{i}(2i\aNd_{A_i}-\X_{A_i}\aDe)\D_{A_{i+1}} \cdots \D_{A_{k-1}}\tU \\
%
&&
\displaystyle
=\big(2i \aDe^{i}\aNd_{A_i}- 2(\sum_{j=1}^i \aDe^{i-j}\aNd_{A_i}\aDe^{j})-\X_{A_{i}}\aDe^{i+1}\big)\D_{A_{i+1}} \cdots \D_{A_{k-1}}\tU \\
%
&&
\displaystyle
=
- \X_{A_{i}}\aDe^{i+1}\D_{A_{i+1}}\cdots \D_{A_{k-1}}\tU .
\end{array}
$$ since $[\aDe,\X_{A_i}]=2\aNd_{A_i}$.  Of course $[\aNd_{A_i},\aDe]$
does not in general vanish. On Ricci-flat manifolds the obstruction is
a sum of terms each of the form $-2\aR_{A_i}{}^P\sharp \aNd_P$ where
the $\sharp$ indicates the usual endomorphism action of the Riemannian
curvature tensor. (Here we view the ambient Riemannian curvature
$\aR_{AB}$ as an ambient 2-form with values in endomorphisms of the
ambient tangent bundle. The abstract indices indicate only the form
indices.) In the further re-expression of terms to the form \nn{form}
and then on to an ambient expression for the curvature terms
polynomial in $\X$, $\aR$, $\D$, $\h$, and $\h^{-1}$ via step 5 of the
algorithm in \cite{GoPetLap} it is easily seen that the free index
$A_i$ remains on an $\aR$. (Note that the only way free indices are
moved from one ambient tensor to another is via the $\aR_{AB}{}\sharp
$ action.)  Via the equivalence of such expressions with tractor
formulae as established in \cite{CapGoamb,GoPetLap} this yields the
claim in the theorem.  \quad $\Box$

In a general scale, the weighted tensor components of $\Pa$ are not
 natural tensors.  Nevertheless it turns out that the $P_{k}$
 operators are often natural, as follows.
\begin{theorem}\label{natthm}
 For a tractor bundle $\ce^\bullet$ 
 the operator $ P_{k}:\ce^\bullet[k-n/2]\to \ce^\bullet[-k-n/2] $ 
is the restriction to conformally
Einstein structures of a natural conformally invariant differential  operator
 if one of
the following holds: $n$ is odd and $k\geq 1 $, or if $n$ even and
$1\leq k\leq n/2-1$, or if $n$ even, $\ce^\bullet[k-n/2]=\ce[k-n/2] $ and
$1\leq k\leq n/2$.
If
$\ce^\bullet[k-n/2]=\ce[k-n/2] $ the operator concerned is a GJMS
operator. Otherwise the operator is a (generalised
GJMS) operator $\Box^0_k$ as in Theorem \ref{powerslap} above.  
\end{theorem}
\noindent{\bf Proof:} Let $\Pa_A=\frac{1}{n}D_A\si$ for and Einstein
scale $\si$. Contracting $\Pa^{A_1}\cdots \Pa^{A_{k-1}}$ into both
sides of \nn{powerform} we obtain
$$
(-1)^{k-1}\si^{n-1} \Box_k^0 U = \Pa^{A_1}\cdots \Pa^{A_{k-1}} \Box D_{A_1}\cdots D_{A_{k-1}}
U .
$$ We have used here that the terms of ${\Cal P}^{\bullet,k}$ are
annihilated. This is clear in view of the result $\Pa^AW_{ABCD}=0$ and
the following two points: in every term of this partial contraction
expressing ${\Cal P}^{\bullet,k}$ at least one of the free indices
$A_1,\cdots ,A_{k-1}$ appears on a $W$; since $\Pa$ is
parallel it commutes with the operators $D$ in the expression for
${\Cal P}^{\bullet,k}$.  Thus in any case where the natural conformally
invariant operators $\Box_k^0$ are defined, as listed in Theorem
\ref{powerslap}, we have that (on conformally Einstein manifolds)
$P_{k}=\Box_k^0$.  \quad $\Box$

\begin{theorem}\label{bigprodthm} 
For an Einstein metric  $g$, $P^g_k:\ce^{\bullet}[k-n/2]\to \ce^{\bullet}[-k-n/2]$ is given by the formula
$$ P^g_{k} = \prod_{l=1}^k(\De^g - b_l J^g),
$$ where $b_l= (n/2+l-1)(n/2-l)(2/n)$, $\Delta^g$ is the
Laplacian for $g$ and $J^g$ the trace of the Schouten tensor for $g$.
\end{theorem} 
\noindent{\bf Proof:} According to the definition \nn{Pdef} we have 
$$ P^g_{k} =(-1)^{k-1}\si^{1-k} \Pa^{A_2}\cdots \Pa^{A_{k}} \Box
D_{A_2}\cdots D_{A_{k}} ,
$$
or equivalently \nn{altform}:
$$
P^{g}_{k} f=(-1)^{k}\si^{-k} \Pa^{A_0}\Pa^{A_1}\cdots
\Pa^{A_{k-1}} D_{A_0} D_{A_1}\cdots D_{A_{k-1}} f .
$$
 Since $\Pa$ is parallel and of weight 0 it commutes with
the tractor D-operator. 
Thus we get 
\begin{equation}\label{powerlapform}
 P_{k}= (-1)^{k}\si^{-k} (\Pa^{A_1}D_{A_1}) \cdots
(\Pa^{A^{k}}D_{A_{k}}) .
\end{equation}

 For an Einstein metric $g$, the parallel tractor is given by
$[\Pa^A]_g=(\si, 0 , -\frac{1}{n}\J \si)$ where $\si\in \ce[1]$ is the
conformal scale corresponding to $g$. That is $g =\si^{-2}\bg$.
Now for $U$ a tractor of weight $k-n/2$, 
$(\Pa^{A_{l+1}}D_{A_{l+1}})\cdots (\Pa^{A_k}D_{A_k}) U $ has conformal weight 
$l-n/2$. If $\tilde{U}$ is a 
tractor weight $l-n/2$
then we have 
$$
\Pa_A D^A \tilde{U} = \si \left(\begin{array}{ccc} -\J/n & 0 & 1 \end{array} \right)  
 \left(\begin{array}{c} (l-1)(2l-n)\tilde{U} \\ 
2(l-1)\nd^a \tilde{U} \\ 
-\Delta \tilde{U} - (l-n/2)\J \tilde{U} \end{array} \right),
$$ 
where on the right-hand-side $\Pa$ and $D^A$ are expressed in terms
of the metric $g$ (i.e. the right-hand-side is $ [\Pa_A]_g [D^A
\tilde{U}]_g$) and $\Delta:=\bg^{ab}\nd_a\nd_b$.  Executing the 
multiplication we obtain
$$
\Pa_A D^A \tilde{U} = -\si (\Delta - (n/2+l-1)(n/2-l)(2/n)\J)\tilde{U}
= -\si (\Delta - b_l\J)\tilde{U} .
$$ 
Once again recalling that $\nd_a \si=0$, from the definition of $\nd$
in the scale $\si$, it follows  that the operator 
$P_{k}:\ce[k-n/2]\to \ce[-k-n/2]$ has the form
$$
P_{k} = \prod_{l=1}^k(\Delta - b_l \J).
$$ This is the result claimed if we view $P_{k} $ as an operator
between density weighted tractors, $ P_{k}:\ce^\bullet[k-n/2]\to
\ce^\bullet[-k-n/2] $.  The equivalent (``covariant'') operator
between weight 0 tractors is given by the composition
$\si^{k+n/2}P_{k}\si^{k-n/2}$ (where we view the powers of $\si$ as
multiplication operators). Once again we can move the $\si$'s around
using that the Levi-Civita connection for $g$ annihilates
$\si$. Writing $\Delta^g =g^{ab}\nd_a\nd_b= \si^2 \Delta$ and
$\si^2\J=J^g=g^{ab}\Rho_{ab}$ the result claimed in the theorem is
immediate.  \quad $\Box$ 

\noindent{\bf Remark:} \label{powerlaprem} From \nn{powerlapform}
above we see that $P_k$ is essentially the ``power'' $
(\Pa^{A_1}D_{A_1}) \cdots (\Pa^{A^{k}}D_{A_{k}}) $ of the
``Laplacian'' $\Pa^AD_A$. The caution to go with this statement is
that $\Pa^A D_A$ has leading term $-\si\Delta$ (rather than $\pm
\Delta$) and, also, since this Laplacian $\Pa^AD_A$ is sensitive to
weight the composition of these Laplacians is not strictly a power.  \quad
\endrk

\noindent{\bf Proof of Theorem \ref{Qthm}:} In \cite{GoPetLap} is
shown that the GJMS operators, i.e.\ the $\Box^0_k$ acting on
densities, factor through the tractor D-operator in the sense that in
\nn{powerform} ${\Cal P}^{\bullet,k}_{A_1\cdots A_{k-1}} U$ (for $U\in
\ce[k-n/2]$) has the form $\tilde{{\Cal P}}^{\bullet,k}_{A_1\cdots
A_{k-1}}{}^B D_B U$. In particular this is true in the case that $n$
is even, $k=n/2$ and $U$ a true function (i.e.\ a density of weight
0). There it is also established that the Q-curvature is then obtained
by replacing, in the right hand-side of this specialisation of
\nn{powerform}, the rightmost $DU$ by the standard tractor $I^g$,
which is given in any scale $g$ by the explicit formula $
I^g_A:=(n-2)Y_A-\J X_A$. (See Section 6.2 of \cite{BrGodeRham} for the
geometric meaning of $I^g$ and its 
relationship to the ambient construction of $Q$ there and in \cite{FeffHir}.)

  That is
$$
\begin{array}{lll}
\lefteqn{\displaystyle(-1)^{n/2-1}X_{A_1}\cdots
X_{A_{n/2-1}} Q}
&&\vspace{2mm}
\\
&&
\displaystyle
=
\Box D_{A_1}\cdots D_{A_{n/2-2}}I^g_{A_{n/2-1}}

+\tilde{{\Cal P}}^{\bullet,n/2}_{A_1\cdots A_{n/2-1}}{}^{B} I^g_{B}
.
\end{array}
$$

In the case that $g$ is Einstein we may contract both sides with the
$(n/2-1)^{\rm st}$-power of the 
corresponding parallel tractor $\Pa$. 
 This yields  
$$
Q=(-1)^{n/2-1}\si^{1-n/2} \Pa^{A_2}\cdots \Pa^{A_{n/2}} \Box
D_{A_2}\cdots D_{A_{n/2-1}}I^g_{A_{n/2}} ,
$$ since from ${\Cal P}^{\bullet,n/2}_{A_1\cdots A_{n/2-1}} f = \tilde{{\Cal
P}}^{\bullet,n/2}_{A_1\cdots A_{n/2-1}}{}^B D_B f$ and Theorem
\ref{powerslap} it is clear that $\Pa^{A_1}\cdots \Pa^{A_{n/2-1}}\tilde{{\Cal
P}}^{\bullet,n/2}_{A_1\cdots A_{n/2-1}}{}^{B} =0$.  
Thus,
once again using that $\Pa$ commutes with the tractor D-operators, we
come to
$$
Q=(-1)^{n/2-1}\si^{1-n/2}  \Box
(\Pa^{A_2}D_{A_2})\cdots (\Pa^{A_{n/2-1}}D_{A_{n/2-1}})(\Pa^{A_{n/2}}I^g_{A_{n/2}}) .
$$ Now $I^g_A$ has conformal weight $-1$ and 
so
$\Pa^{A_{n/2}}I^g_{A_n/2}= 2(1-n)\si \J/n$, where $\J$ is the $\bg$-trace of the Schouten tensor for $g$. Arguing otherwise as in the
proof of the previous theorem gives
$$
Q=\frac{2(1-n)}{n}\Big(\prod_{l=1}^{n/2-1}(\Delta - b_l \J)\Big)\J ~.
$$
Since $\nd \J=0$ this simplifies at once to 
$$
Q=(-1)^{n/2}\frac{2(n-1)}{n}\J^{n/2}\prod_{l=1}^{n/2-1}b_l .
$$ Here we have viewed $Q$ as a density of weight $-n$. But
multiplying by $\si^n$ and then absorbing these via $J^g=\si^2\J=g^{ab}\Rho_{ab}$ we
obtain the same formal expression except we will write $Q^g$ to mean the
 function equivalent to $Q$.  Finally using that $\J^g={\rm Sc}^g/(2(n-1))$,
inserting $b_l= (n/2+l-1)(n-2)/n$ and re-arranging we obtain
\nn{Qformula}.  \quad $\Box$

\section{Naturality and related issues ``beyond the obstruction''} 
\label{farside}

Theorem \ref{natthm} established that in odd dimensions, and for
sufficiently small $k$ in even dimensions, the $P_{k}$ are natural
conformally invariant differential operators. On the other hand
Theorem \ref{bigprodthm} shows that when expressed in terms of an
Einstein metric $g$ then the $P^g_{k}$ are given by natural (and very
simple) formulae, for all positive integers $k$. The remaining issue
is whether, in the case of even dimensions and ``large $k$'' (i.e. $k$
outside the range covered by Theorem \ref{natthm}), the $P_k$ may be
given by a natural formula in terms of an arbitrary metric from the
(conformally Einstein) conformal equivalence class. Note that on
conformally flat structures there is such a formula in all cases.

Here we consider the naturality issue and an
explicit verification of Theorem \ref{g1g2} for a specific case viz.\
$P_{6}$ acting on densities. Note that in dimension 4 this is a large
$k$ (in the sense described above) case. Let us first consider all other
dimensions. From \cite{GoPetLap} we have the formula
\begin{equation}\label{6th}
(n-4)X_AX_B \Box^0_3 f = (n-4)\Box D_A D_B f + 2 W_{A}{}^C{}_B{}^E D_C D_E f 
\end{equation}
for $f\in \ce[3-n/2]$ and dimensions $n\neq 4$. 
This formula holds for any conformal structure.
Now suppose that we have a conformally Einstein manifold and that  $\Pa$ is an
Einstein tractor corresponding to a(n Einstein) scale $\si$. 
Contracting $\Pa$ into both free tractor indices of \nn{6th} gives
$$
(n-4) \Box^0_3 f = (n-4)\si^{-2}\Pa^A\Pa^B \Box D_A D_B f =(n-4)P_3 f~,
$$
since (from \nn{PaW})
\begin{equation}\label{IIW}
\Pa^A\Pa^B W_{A}{}^C{}_B{}^E =0,
\end{equation}
and $\Pa^AX_A=\si$.
Thus in dimensions other than 4 we have $\Box^0_3 f = P_3 f$
 and $P_3 f$ is a 
natural operator as claimed in Theorem \ref{natthm}.

To simplify the subsequent 
discussion let us calculate in some scale $g$ and agree to write all
natural operators as contractions or partial contractions involving
universal polynomials in (Levi-Civita) covariant derivatives with
coefficients depending polynomially on the conformal metric, its
inverse, the Weyl tensor tensor, the Schouten tensor and their
covariant derivatives. With these conventions note that, from the
formulae \nn{Dform} and \nn{Wform}, the tractor D-operator and the
W-Tractor are polynomial in the dimension $n$.
Thus the right hand of \nn{6th}
side is given by a universal formula  that is polynomial in $n$. 

First let us observe it is straightforward to extract from \nn{6th} a
formula for $\Box^0_3$ (in dimensions $n\neq 4$). This formula sheds
some light on the difficulties in dimension 4. Using \nn{Wform},
\nn{Dform} and \nn{basictrf} we have
\begin{equation}\label{invt}
Y^BY^C(2 W_{B}{}^S{}_C{}^T D_S D_T f )=8B^{cd}(2\nd_c\nd_d f-(n-6)\Rho_{cd}f).
\end{equation}
So from \nn{6th} we obtain the expression
\begin{equation}\label{Neq}
 \Box^0_3 f =  N f + \frac{8}{n-4}B^{cd}(2\nd_c\nd_d f-(n-6)\Rho_{cd}f),
\end{equation}
where we note that the differential operator $N f:=Y^AY^B \Box D_A D_B
f $ is polynomial in $n$ and natural. (For the last point note that
from \nn{Dform} it is clear that $D$ is natural therefore so are its
iterations. Contracting with $ Y^AY^B$ simply extracts a component
which is necessarily natural by the definition of naturality for
tractor operators.)  In general, the term on far right of the display
is undefined in dimension 4, and so, in some sense, it is this that
prevents $\Box^0_3 $ from yielding a conformally invariant operator in
dimension 4.  It follows from \nn{Wform} that in dimension 4 the
expression \nn{invt} is a (coupled) conformal invariant.

Observe that in the expression \nn{Wform}, for $W_{ABCE}$, each term
has a factor of $(n-4)$ except for one which involves the Bach tensor
$B_{eb}$.  Let us suppose now, and henceforth, that we are on a
conformally Einstein structure and that $\Pa$ is a parallel tractor
corresponding to an Einstein metric $g=\si^{-2}\bg$. In this case we
have \nn{PaW}, and expanding the identity $W_{BCDE}\Pa^E =0$ (see
\cite{GoNur} for details) we obtain that 
$$
B_{cd}=(4-n)\si^{-1} A_{cde}\nd^e\si .
$$ 
Substituting this into \nn{Wform} we obtain that 
$$
W_{ABCE}=(n-4)\tilde{W}_{ABCE}
$$
where $\tilde{W}_{ABCE}$ is defined to be
$$
\begin{array}{l}
 Z_A{}^aZ_B{}^bZ_C{}^cZ_E{}^e C_{abce}
-2 Z_A{}^aZ_B{}^bX_{[C}Z_{E]}{}^e A_{eab} \\ 
-2 X_{[A}Z_{B]}{}^b Z_C{}^cZ_E{}^e A_{bce} 
- 4 X_{[A}Z_{B]}{}^b X_{[C} Z_{E]}{}^e \si^{-1} A_{eba}\nd^a\si~.
\end{array}
$$ 
Note that we have $\Pa^A\Pa^B \tilde{W}_{A}{}^C{}_B{}^E =0$ as a
tautological identity.  
Now consider
\begin{equation}\label{divver}
\Box D_A D_B f + 2 \tilde{W}_{A}{}^C{}_B{}^E D_C D_E f .
\end{equation}
It follows from \nn{6th} that, in dimensions $n\neq 4$, this is an
expression for $X_AX_B \Box^0_3 f$.  Thus, in these dimensions, it has
the property that contraction by a product of any pair of projectors
($X^AX^B,$ $X^AZ^B{}_b$ and so forth) annihilates it, with the
exception of $Y^AY^B$. Now as an operator on $f$ {\em and} $\si$,
\nn{divver} is a natural operator given by a universal formula
polynomial in the dimension $n$. It follows that, in fact, in all
dimensions it has the property that contraction by a product of any
pair of projectors, other than $Y^AY^B$, annihilates it.  Thus
\nn{divver} defines an operator $\tilde{P}_3$ on $\ce[3-n/2]$ by
$$
X_AX_B \tilde{P}_3 f= \Box D_A D_B f + 2 \tilde{W}_{A}{}^C{}_B{}^E D_C D_E f .
$$ 
Now contracting $\Pa$ onto both free indices and using that
$\Pa^A\Pa^B \tilde{W}_{A}{}^C{}_B{}^E =0$ we see that in all
dimensions $\tilde{P}_3=P_3$. On the other hand, using this and contracting the
projector $Y$ onto both free indices we obtain an explicit formula for $P_3$
\begin{equation}\label{P6ex}
 P_3 f= N f - 8   A^{cde}(\si^{-1}\nd_e\si)  (2\nd_c\nd_d f-(n-6)\Rho_{cd}f) ,
\end{equation} 
where, as above, $N f=Y^AY^B \Box D_A D_B f $.  The right-hand side
of the display is given by a natural expression if we calculate in the
Einstein scale $\si$ (i.e.\ $\si^{-2}\bg$ so that $\nd_e\si=0$) but
apparently, in dimension 4, cannot be re-expressed as a natural
formula for a general metric $g$ from the (conformally Einstein)
conformal class.

Finally note that if $(\si_1)^{-2}\bg$ and $(\si_2)^{-2}\bg$ are two Einstein 
metrics then 
$$
A^{cde}(\si_2\nd_e\si_1-\si_1\nd_e\si_2)=0,
$$ where the Cotton tensor $A_{cde}$ is for any metric $g$ in the
conformal class. Using \nn{tractcurv} this is an easy consequence of
$(\si_2\nd^a\si_1-\si_1\nd^a\si_2) \Om_{abCD}=0$, from Proposition
\ref{kdotCetc}.  Thus, in \nn{P6ex}, $A^{cde}(\si^{-1}\nd_e\si)$ is
independent of which Einstein scale we use. (Note that in dimensions
other than 4 this is already clear from the identity
$B_{cd}=(4-n)\si^{-1} A_{cde}\nd^e\si$ since $B_{cd}$ is natural.)
This verifies explicitly in this case the general result of Theorem
\ref{g1g2}.

\end{document}